\documentclass[11pt]{article}
\usepackage{amsmath,amssymb, amsfonts, amsthm, latexsym, hyperref,titlesec }
\usepackage{bbm,xcolor,graphicx}

\titleformat{\subsubsection}[runin] {\normalfont\bfseries}{\thesubsubsection}{0.7em}{\addperiod}
\newcommand{\addperiod}[1]{#1.}

\titlespacing\subsubsection{5pt}{4pt plus 4pt minus 2pt}{5pt plus 2pt minus 2pt}

\newtheorem{question}{Question}

\newcommand{\N}{\mathbb{N}}
\newcommand{\Z}{\mathbb{Z}}

\begin{document}
\title{$p_c$, $p_u$ and graph limits }
\author{Itai Benjamini }

\date{2.7.17}

\maketitle

\begin{abstract}
A question relating the critical probability for percolation, the critical probability for a unique infinite cluster
and graph limits is presented, together with some partial results.
\end{abstract}

\footnotetext{MSC subversive math}

\section{A question }

Let $G$ be an infinite Cayley graph rooted at $o$.
Denote by $B(o, n)$ the ball of radius $n$ centered at $o$, and let $H$
be the limit, in the sense of \cite{BS},
of the sequence of balls $B(o, 1),B(o,2)...$,
possibly passing to a subsequence to get convergence.
$H$ is a unimodular random graph, \cite{AL}.

\begin{question}
\label{c1}
Does $p_c(H) = p_u(G)  \mbox{   } a.s.$ ?
\end{question}

In $p$-Bernoulli site percolation on a graph, each vertex is open with probability $p$ independently.
$p_c$ denotes the critical probability for the existence of an infinite open cluster (connected component),
and $p_u$ is  the critical probability for uniqueness of the infinite open cluster,  see \cite{BS96}.
\medskip

A limit, a la \cite{BS}, of finite graphs $G_n$  is a random rooted infinite graph
with the property that neighborhoods of $G_n$ {\it around a random uniform vertex} converge
in distribution to neighborhoods of the infinite graph around the root. When $G$ is a Cayley graph
it might be the case that the limit along balls exists, and there is no need to pass to a subsequence?
\medskip

Neither one of the directions in question~\ref{c1} seems immediate.
For a very simple example take $G$ to be the $3$-regular tree, $T_3$, then $H$ is the canopy tree,
which can be thought of as an ``infinite binary tree viewed from a leaf".
$p_c(H) = p_u(G) =1$. While for amenable Cayley graphs $p_c = p_u$ \cite{BK}, it is conjectured that $p_c < p_u$
for non-amenable Cayley graphs \cite{BS96}. It was shown in \cite{PS} that every non-amenable group admits a generating set
for which $p_c < p_u$. Showing that $p_c(H) > p_c(G)$ for  non-amenable $G$, will be  an interesting progress towards question~\ref{c1}.
In \cite{LS}  it is proved that $p > p_u$ is equivalent to non-decay of the connection probability.
For an amenable Cayley graph $G$, the limit along a F{\o}lner sequence is a.s. $G$.
\medskip

\medskip

Recall the critical parameters for weak survival of the contact process
$\lambda_c$ and of strong survival $\lambda_s$, see \cite{P} \cite{L}.

\begin{question}
\label{c2}
Does $\lambda_c(H) = \lambda_s(G)  \mbox{   }a.s.$ ?
\end{question}

Question~\ref{c2}  is already non trivial for $T_d$, and it will be great to see a proof.
\medskip

\section{Partial results and further questions}

Given a Cayley graph $G$, it is often hard to understand the structure of $H$.
Still we now present several families of graphs for which the answer to question~\ref{c1} is yes.
\medskip

\noindent
{\it Planar triangular lattices}

When $G$ is a cocompact triangular lattice in the hyperbolic plane, $H$ is a.s. a triangulated quasi-canopy tree.
That is, a tree quasi isometric to the canopy tree. For simplicity we will assume it's the canopy tree.
View the canopy tree as an infinite ray with growing trees rooted along the ray. For every $p_c(G) < p < p_u(G)$,
and any fixed vertex $v$  on the ray, there is a positive probability  of an open path from $v$ to the leaves.
This gives a cutset separating $o$ from infinity. When $p > p_u(G)$,  the probabilities
of $(1-p)-$ open paths from $v$ to the leaves decays exponentially with the distance from $v$ to $o$.
As the open path has to exit a ball of radius $dist(v, o)$ centered at $v$,
which is isometric to a ball in the transitive lattice $G$.
Therefore restricted to the ball this  event has identical probability to the event of
an open path of a subcritical percolation  on $G$, \cite{BS01} exiting a ball of radius $dist(v, o)$, and by
\cite{DT} this probability decays exponentially.
\medskip

Verify question~\ref{c1}  for  planar stochastic hyperbolic triangulations, \cite{C}?
\medskip

\noindent
{\it Graph products}

It is interesting  to consider next the conjecture for $T_d  \times  \Z$, having a non planar instance with $p_u < 1$ analyzed.
Tom Hutchcroft challenged us with the following comment.
Grimmett and Newman's paper \cite{GN} implies the following statement:
Take the product of the $d$-regular tree with $\Z$, and replace every $\Z$-edge with a path of length $n$.
Then no matter how large $n$ is, $p_u$ of the resulting graph is always bounded from above by $(d-1)^{-1/2}$.
If question~\ref{c1}  has a positive answer, that would imply that if we take the product of the canopy tree with $\Z$
and stretch the lengths of the $\Z$-edges  as much as we want, then $p_c$ is always bounded above by $(d-1)^{-1/2}$ also.
Tom (private communication) proved that this is indeed the case. This supports a positive answer for $T_d \times \Z$.
A first step towards a positive answer to question~\ref{c1} suggested by Tom, is to show that
$p_u(T_d \times \Z) = p_u(T_d \times \N)$. Is there non-uniqueness at $p_u$ on $T_d \times \N$?
\medskip

Percolation on $T_d \times \Z$ at $p_u$ admits infinitely many infinite clusters a.s., while on planar hyperbolic lattices
there is a unique infinite cluster at $p_u$, \cite{Sh}, \cite{BS01}. What is the behaviour at $p_c$ for the corresponding $H$'s?

Does $p_c(\mbox{canopy} \times \Z) = p_c(\mbox{canopy} \times \N)$? Show no percolation at criticality on $\mbox{canopy} \times \N$?

\medskip

When a Cayley graph $G$ is transient for the simple random walk, one expects that the infinite clusters of Bernoulli percolation are a.s.
transient as well. This is known for many families of Cayley graphs but not in general.
Show that the infinite cluster of supercritical Bernoulli  percolation on $ \mbox{canopy} \times  \Z$ is a.s. transient?
\medskip

\noindent
{\it Free products}

In \cite{BE} it was suggested that for any infinite Cayley graph $G$,
the graph limit along a converging subsequence of balls, $H$, is invariantly amenable,
and therefore $p_c(H) = p_u(H)$ and $\lambda_c(H) = \lambda_s(H)$.
In \cite{S} it is proved that if $G$ has Yu's property A, then $H$ is hyperfinite or invariantly amenable.
These include hyperbolic groups, groups with finite  asymptotic dimension and amenable groups.
(\cite{BE} contains many further questions regarding graph limits along balls).
\medskip

Other Cayley graphs Tom Hutchcroft  suggested to look at, as a possible counterexamples, are free product of $\Z^2$ and an edge,
and  a free product of a hyperbolic Cayley graph and an edge.

The free product of $\Z^2$ and an edge has finite asymptotic dimension, \cite{BDK}.
Hence the limit along balls is invariantly amenable
(or hyperfinite) unimodular random graph  and thus the limiting graph has one end a.s., see \cite{AL}.
$H$ is not contained in a finite number of copies of $\Z^2$  as
$H$ has a uniform exponential lower bound on it's volume growth, as $G$ contained a regular tree and therefore
$H$ contained a canopy.

$H$ intersects any copy of $\Z^2$ in a finite set, as otherwise we will get more than one end.
Thus $p_c(H) =1$ a.s., since the infinite  cluster has to cross infinitely many cut edges.

The free product of an hyperbolic group $G$ and an edge is hyperbolic.
Limits along balls of hyperbolic groups is invariantly amenable,
So the limit has one end a.s.
and hence  $p_c(H) = 1$ a.s., by the same argument as above.

\begin{question}
\label{q1}
Assume a Cayley graph $G$ has infinitely many ends. Show that $p_c(H) =1$.
\end{question}

To answer question~\ref{q1}, one needs to show that the limit along balls of any Cayley graph
which is  not quasi isometric to $\Z$, has one end a.s.
\medskip

\noindent
{\it Final comments}

It is still open if there is a Cayley graph in which the sequences  of balls
$B(o, 1),B(o,2)...$ is an expander family, \cite{B}? I conjecture that there is no such Cayley graph.
In this case $H$ will be non-amenable.

Recall a sequence of finite subgraphs  in a graph  is exhausting,
if  each graph in the sequences contains all the previous subgraphs and the union
is the whole graph.

If $G$ is hyperbolic and we consider the limit along an arbitrarily
sequence of  exhausting finite sets, we expect to get $p_c(H) \geq p_u(G)$.
This might hold for any  exhausting sequence of   finite sets
on any Cayley graph. Gromov's monster, \cite{G}, might serve as a counter example.

\medskip

Let $G$ be the Diestel-Leader graph $DL(3,2)$, \cite{DL},  which is a non-amenable, non-unimodular  vertex transitive graph.
The limit along balls of  $DL(3,2)$  is a horocyclic product of two
canopy trees, see figure. Show that $p_c(H) < 1$. Does $p_c(H) = p_u(G)$?
Does $p_c(G) < p_u(G)$?
\medskip

\begin{figure}[h]
\begin{minipage}[b]{1\linewidth}
\centering
\includegraphics[scale=0.14]{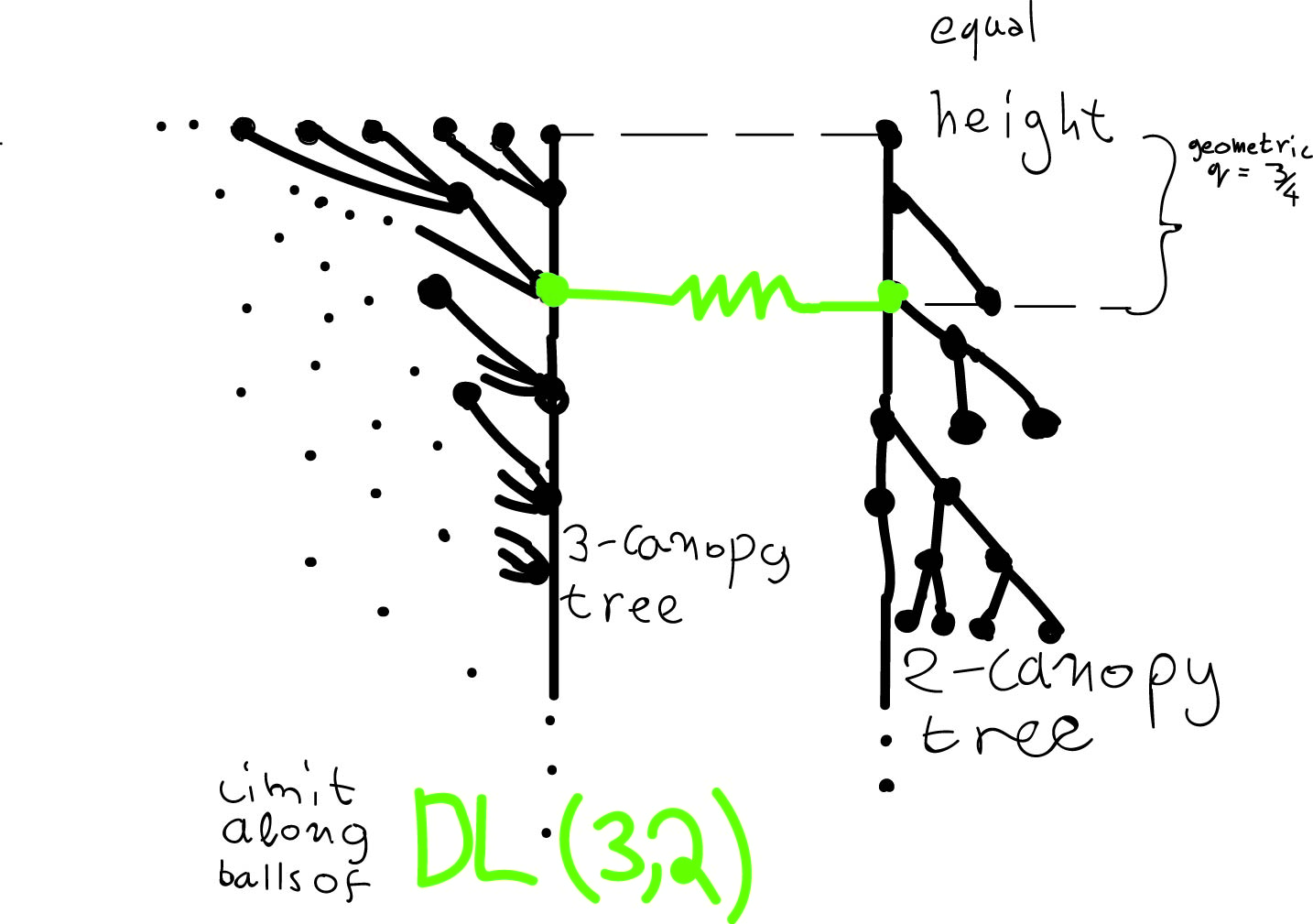}
\label{fig:DL5}
\caption{ The graph limit along balls of $DL(3,2)$.}
\end{minipage}
\quad
\end{figure}

\begin{figure}[h]
\begin{minipage}[b]{1\linewidth}
\centering
\includegraphics[scale=0.14]{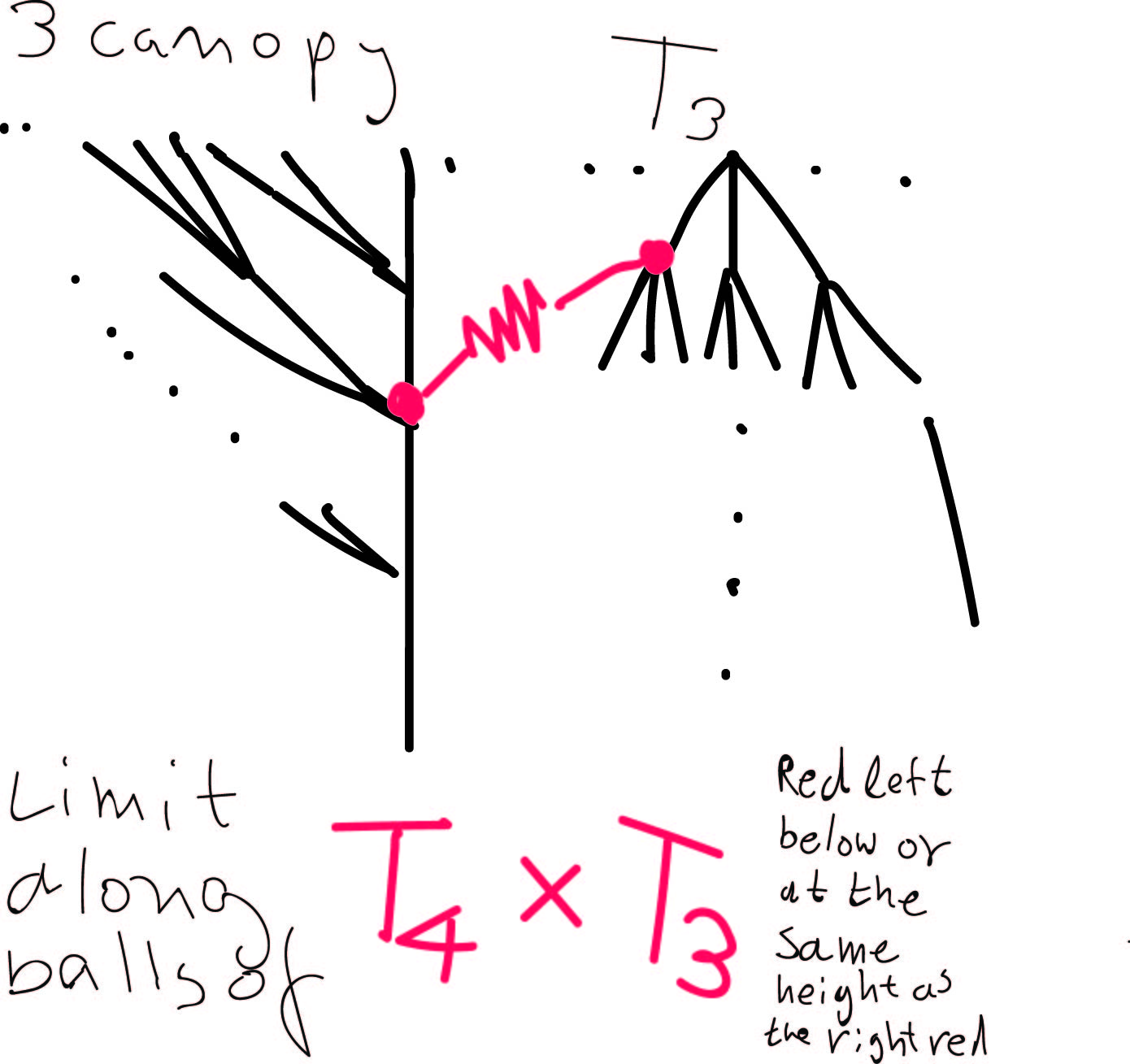}
\label{fig: }
\caption{ The graph limit along balls $T_4 \times T_3$.}
\end{minipage}
\quad
\end{figure}

Maybe a version of question~\ref{c1} for the critical temperature for  magnetization of the ferromagnetic Ising model,
will be more tractable?
\medskip

We believe that the answer to question~\ref{c1}  will still be positive   if rather than taking limit along balls,
we  consider "round" sets such as limit of isoperimetric minimizers, or minimizers along Green balls,
or taking a sequence sets $\{K_n\}$, where $K_n$ is a set of vertices of size $n$ minimizing
the expected escape probability of simple random walk starting from a random uniform element on the set.

Assume $G$ and $G'$ are quasi isometric. Let $H$ be the limit along some exhausting sequence of $G$.
Is there an exhausting sequences of $G'$, so that $H'$ is quasi isometric to $H$? That is, there is a coupling of the two
random rooted graphs $H$ and $H'$ which is a quasi isometry.

\medskip

We would like to end with a random walk question of similar spirit.
In \cite{BS} it was proved that the limit of bounded degree finite planar graphs, is a.s. recurrent
for the simple random walk.

Which Cayley graphs has an exhausting sequence of subgraphs with a limit that is a.s. transient?
One feature of transient bounded degree  planar graphs, is that they admit non constant Dirichlet harmonic functions,
that is, harmonic functions with a gradient in $l^2$,
see \cite{BS96} and \cite{H}. Let $G$ be an infinite transient Cayley graph.

\begin{question}
\label{c4}
Assume $G$ has no non constant Dirichlet harmonic functions. Is there an exhausting  sequence of finite subgraphs, such that
$H$ is a.s. transient?
\end{question}

If there is such a sequence,  probably the limit along growing balls is a.s. transient.

There are no non constant Dirichlet harmonic functions on $T_3 \times \Z$.
The limit along boxes of $T_3 \times \Z$ is a.s. a $\mbox{canopy tree} \times \Z$, which is transient.
\medskip

All examples we know, including limits along balls in  products of trees, are Liouville,
 
\begin{question}
Are there examples in which limits along balls are a.s. non-Liouville?
\end{question}
\medskip

\noindent
{\bf Acknowledgements:} Thanks to Lewis Bowen,  Tom Hutchcroft for  very useful discussions.

\end{document}